\documentclass[a4paper,10pt,twoside]{article}

\usepackage[body={165mm,235mm}]{geometry}
\usepackage{a1}
\usepackage[english]{babel}
\usepackage[latin1]{inputenc} 
\usepackage{amsfonts,amssymb}
\usepackage{amsmath}
\usepackage{graphicx}	

\usepackage{makeidx}
\makeindex

\def\mapright#1#2#3{\smash{\mathop{\hbox to
#3{\rightarrowfill}}\limits^{#1}_{#2}}}

\def\mapleft#1#2#3{\smash{\mathop{\hbox to
#3{\leftarrowfill}}\limits^{#1}_{#2}}}

\def\mapright#1#2{\smash{\mathop{\hbox to 0.90cm{\rightarrowfill}}\limits^{#1}_{#2}}}
\def\mapleft#1#2{\smash{\mathop{\hbox to 0.90cm{\leftarrowfill}}\limits^{#1}_{#2}}}

\def\mapleftright#1#2{\smash{\mathop{\hbox to 0.80cm{\leftarrowfill \rightarrowfill}}\limits^{#1}_{#2}}}

\title{A challenge to 3-manifold topologists and group algebraists
} 
\author{Sóstenes L. Lins and Lauro D. Lins}

\date{\today}

\begin{document}

\maketitle

\begin{abstract}
This paper poses some basic questions about instances (hard to find) of 
a special problem in 3-manifold topology.
``Important though the general concepts and propositions may be with the 
modern industrious passion
for axiomatizing and generalizing has presented us \ldots nevertheless I am convinced that
the special problems in all their complexity constitute the stock and the core of mathematics; 
and to master their difficulty requires on the whole the harder labor.'' Hermann Weyl 1885-1955, 
cited in the preface of the first edition (1939) of \cite{whitehead1997}.
\end{abstract}

\section{A doubt in the classification of 3-manifolds: $U[1466]$ and $U[1563]$}

The objective of this note is to pinpoint an aspect of the classification of 
3-manifolds which is very important and has been essentially neglected in
the last 35 years of successes with the work of W. Thurston, G. Perelman,
I. Agol and many others. In despite of enormous progress, the 
classification problem remains, to our eyes, very difficult. The aspect we want to pinpoint is
asking basic questions on hard to find tough instances of the general theory.

\begin{figure}[!h]
\begin{center}
\includegraphics[width=10cm]{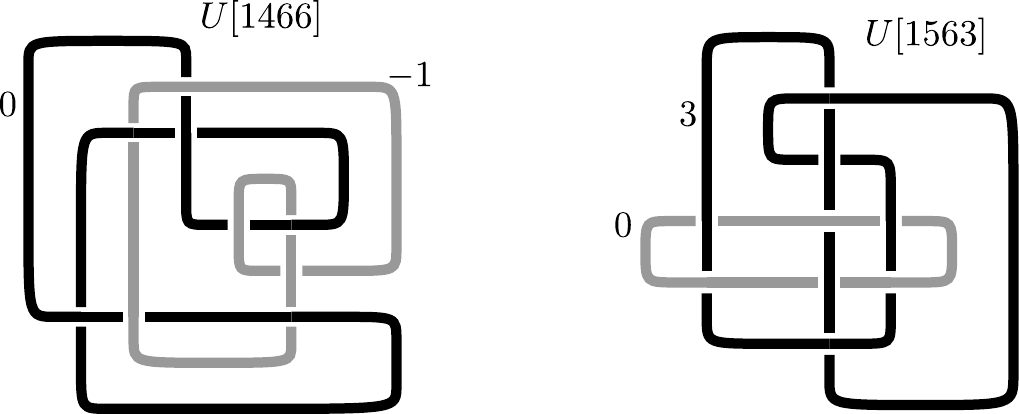}
\caption{\sf Are the closed orientable 3-manifolds obtained from surgery on $\mathbb{S}^3$ of the above blackboard 
framed links followed by the canonical Dehn fillings homeomorphic, or not?}
\label{fig:firstdoubtA}
\end{center}
\end{figure}

In a fundamental paper W. B. R. Lickorish proves that each closed orientable 
3-manifold can be encoded as a link in $\mathbb{S}^ 3$ with integers in 1-1 correspondence
with its components, \cite{lickorish1962representation}, the so called {\em framed links}.

Consider the two closed orientable 3-manifolds
obtained from surgery and canonical Dehn fillings on the 2-component
blackboard framed \cite{kauffman1991knots} link of Fig. \ref{fig:firstdoubtA}. 
Both are homology spheres, so their fundamental groups are perfect.
SnapPea \cite{weeks2001snappea} tells us, according to S. Matveev \cite{matveev2008},
that they are both hyperbolic and 
have the same volume up to 10 decimal places,. Moreover,
their Witten-Reshetiken-Turaev invariants with 10 decimal places 
agree up to $r=12$, \cite{lins2007blink}. These facts seem to imply that the manifolds are homeomorphic.
However, computations based on the methodology of 
\cite{lins2007blink} and \cite{lins1995gca}, which were up to this point successful
in finding homeomorphism between pairs of 3-manifolds, appear to fail for the first time.
Our bet is that the methodology does not fail, that is, the manifolds are not homeomorphic.
In the last 5 years we have asked the help of 
various distinguished topologists in trying to settle
this example. None of them succeeded in answering our question. So,
we believe the time is ripe to bring our doubt to the broader community 
of mathematicians
dealing with 3-manifolds and/or combinatorial group theory. 
This example corresponds to
the pair of blackboard framed links $U[1466]$ and $U[1563]$ 
of \cite{lins2007blink}. The numbers attached to the components (framing) 
coincide with their self-writhes in the given projection and, so, can be discarded.
Note that by introducing an appropriate number of positive or negative curls we can obtain
any framed link as a blackboard framed link (and discard the framings). In a blackboard framed
link we do not need nor use the framing to obtain a presentation of the fundamental group.

If the manifolds being compared are hyperbolic, then 
the difficult topological question of homeomorphism between the manifolds transforms
into the possibly equally difficult algebraic question of
isomorphism between their fundamental groups. So, as long as the general associated
question is not settled,
we have replaced a problem which we do not know how to solve into another,
which we also do not know how to solve. This might be, in some aspects,
progress, but hardly a definitive one.
In general, how to prove that the fundamental groups of hyperbolic 3-manifolds
are not isomorphic? Start by proving that there is no isomorphism between the 
fundamental groups of the above 3-manifolds. Or find one. 



\begin{figure}[!h]
\begin{center}
\includegraphics[width=12.5cm]{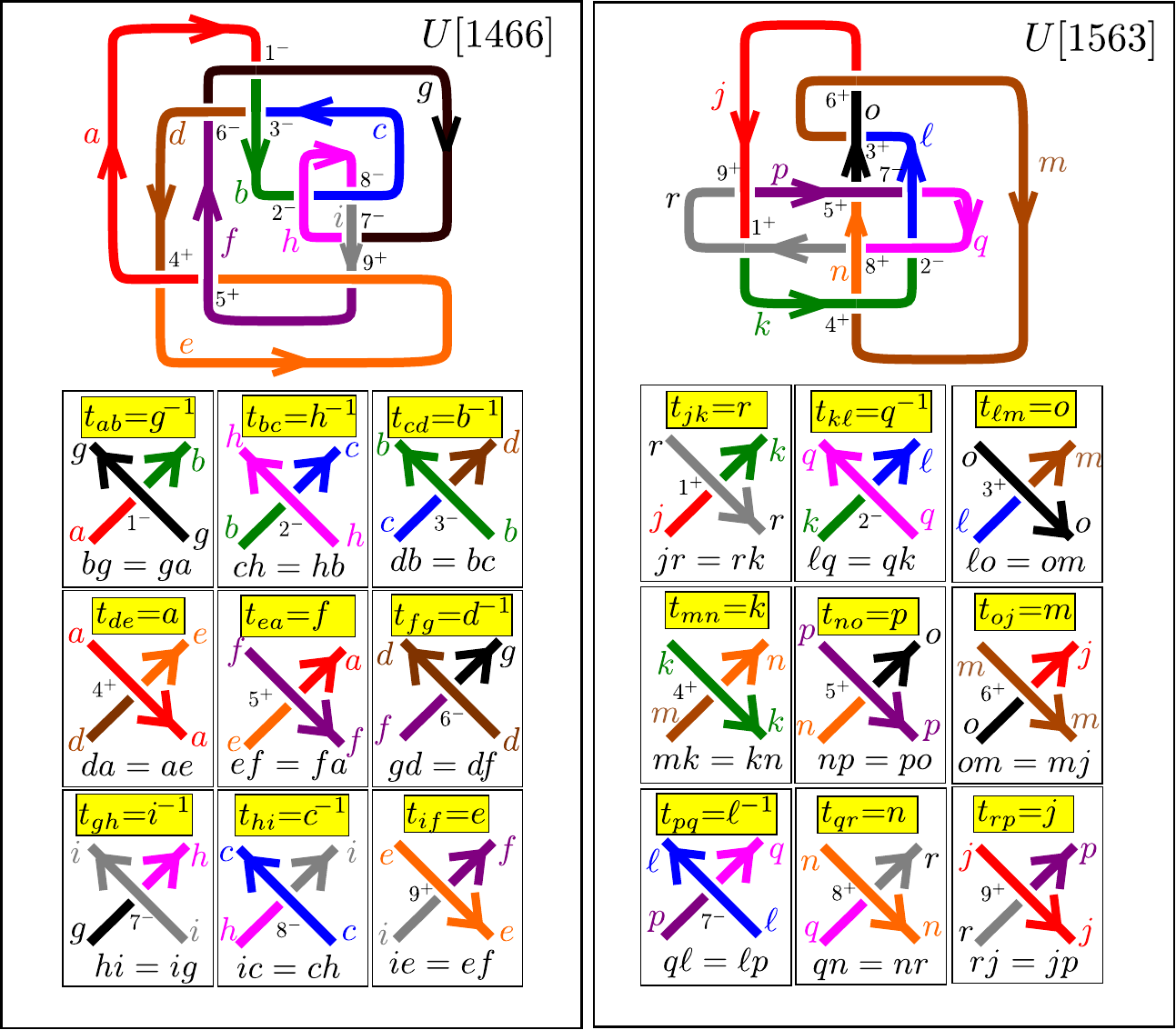}
\caption{\sf Finding presentations for the fundamental groups  of $M^ 3[1466]$ and $M^ 3[1563]$:
we arbitrarily orient the links, write the transition generators, $t_{xy}$'s,  in terms of the Wirtinger 
generators (\cite{stillwell1993classical}), write the Dehn fillings relators (\cite{rolfsen2003knots}) in 
terms of the transition generators and, finally, write the Wirtinger 
relations for the fundamental groups of the exterior of the links.}
\label{fig:firstdoubtB}
\end{center}
\end{figure}

The presentations for the fundamental groups of the manifolds $M^3[1466]$ and $M^3[1563]$ are:
\begin{center}
$\pi_1[1466]=\langle \{t_{ab},t_{bc},t_{cd},t_{de},t_{ea},t_{fg},t_{gh},t_{hi},t_{if},a,b,c,d,e,f,g,h,i\},$\\
$\{t_{ab}=g^{-1},t_{bc}=h^{-1},t_{cd}=b^{-1},t_{de}=a,t_{ea}=f, $\\
$t_{fg}=d^{-1},t_{gh}=i^{-1},t_{hi}=c^{-1},t_{if}=e,$\\
$t_{ab}t_{bc}t_{cd}t_{de}t_{ea}=1 ,t_{fg}t_{gh}t_{hi}t_{if}=1,$\\
$bg=ga, ch=hb, db=bc,da=ae, ef=fa, gd=df, hi=ig,ic=ch,ie=ef\} \rangle,$
\end{center}
\begin{center}
$\pi_1[1563]=\langle \{t_{jk},t_{kl},t_{lm},t_{mn},t_{no},t_{oj},t_{pq},t_{qr},t_{rp},j,k,l,m,n,o,p,q,r\}$\\
$\{t_{jk}=r,t_{kl}=q^{-1},t_{lm}=o,t_{mn}=k,t_{no}=p,t_{oj}=m, $\\
$t_{pq}=l^{-1},t_{qr}=n,t_{rp}=j,$\\
$t_{jk}t_{kl}t_{lm}t_{mn}t_{no}t_{oj}=1,t_{pq}t_{qr}t_{rp}=1,$\\
$jr=rk,lq=qk,lo=om,mk=kn,np=po,om=mj,ql=lp,qn=nr,rj=jp\} \rangle.$
\end{center}

\section{Another doubt: $U[2125]$ and $U[2165]$}

\begin{figure}[!h]
\begin{center}
\includegraphics[width=10cm]{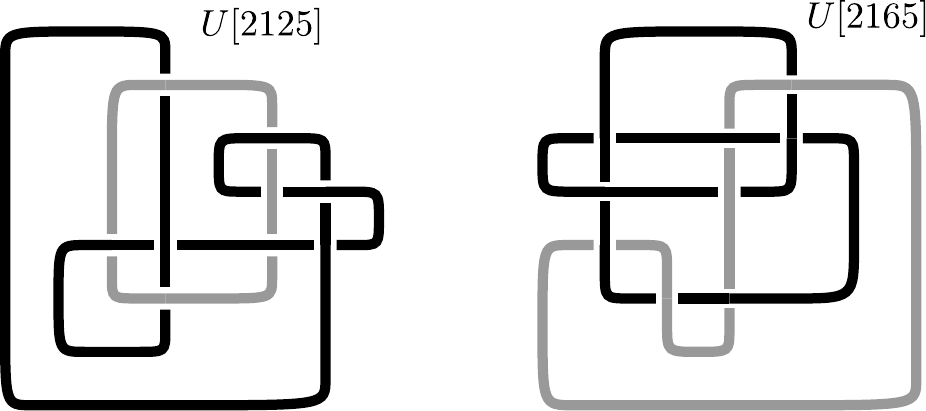}
\caption{\sf Are the closed orientable 3-manifolds obtained from surgery on $\mathbb{S}^3$ of the above
blackboard framed links followed by canonical Dehn fillings homeomorphic, or not? 
The framing of a component in the above links 
is its self-writhe in the given projection. }
\label{fig:seconddoubtA}
\end{center}
\end{figure}

It is important also to distinguish the pair 3-manifolds induced by the blackboard framed links of
Fig. \ref{fig:seconddoubtA}. 
As the previous pair, they are closed hyperbolic homology spheres and their WRT-invariants
agree up to $r=12$ with 10 decimal places, \cite{lins2007blink}.
The presentations for the fundamental groups of the manifolds $M^3[2125]$ and $M^3[2165]$
are:\\

\begin{center}
$\pi_1[2125]=\langle \{t_{ab}, t_{bc}, t_{cd}, t_{de},t_{ef}, t_{fa}, t_{gh}, t_{hi}, t_{ig},a,b,c,d,e,f,g,h,i\},$\\
$\{t_{ab}=h^{-1}, t_{bc}=d, t_{cd}=g^{-1}, t_{de}=b, t_{ef}=a, t_{fa}=i,$\\
$t_{gh}=c^{-1}, t_{hi}=f, t_{ig}=e^{-1},$\\ 
$t_{ab}t_{bc}t_{cd}t_{de}t_{ef}t_{fa}=1, t_{gh}t_{hi}t_{ig}=1,$\\
$bh=ha,bd=dc,dg=gc,db=be,ea=af,fi=ia,hc=cg,hf=fi,ge=ei\} \rangle,$
\end{center}

\begin{center}
$\pi_1[2165]=\langle \{t_{jk},t_{kl},t_{lm},t_{mn},t_{no},t_{oj},t_{pq},t_{qr},t_{rp},j,k,l,m,n,o,p,q,r\},$\\
$\{t_{jk}=r^{-1},t_{kl}=q,t_{lm}=j^{-1},t_{mn}=k,t_{no}=p^{-1},t_{oj}=l^{-1}, $\\
$t_{pq}=n^{-1},t_{qr}=m,t_{rp}=o^{-1}, $\\
$t_{jk}t_{kl}t_{lm}t_{mn}t_{no}t_{oj}=1,t_{pq}t_{qr}t_{rp}=1,$\\
$\{kr=rj,kq=ql,mj=jl,mk=kn,op=pn,jl=lo,qn=np,qm=mr,po=or\} \rangle.$
\end{center}

\begin{figure}[!h]
\begin{center}
\includegraphics[width=12.5cm]{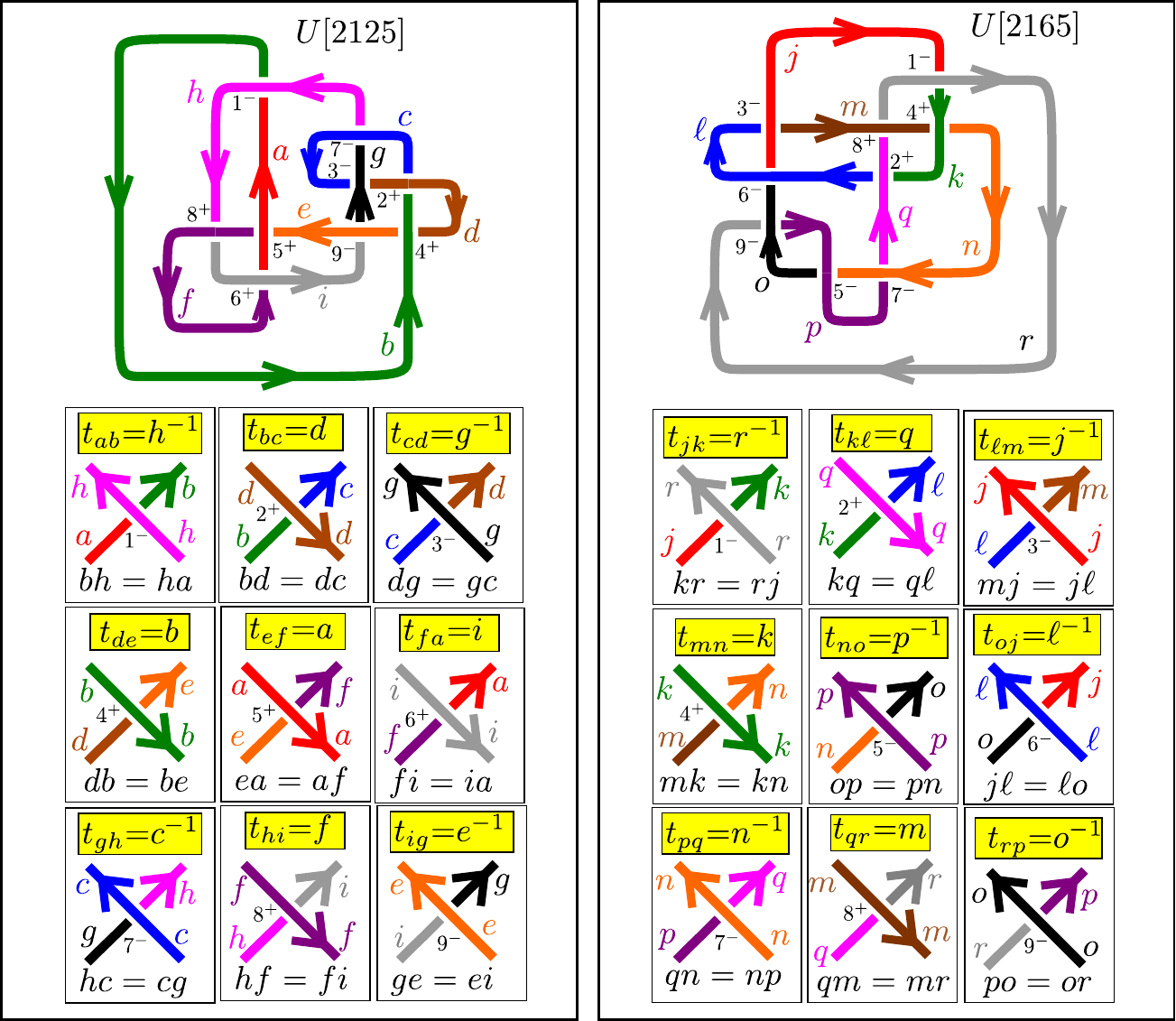}
\caption{\sf Finding presentations for the fundamental groups of $M^ 3[2125]$ and $M^ 3[2165]$}
\label{fig:seconddoubtB}
\end{center}
\end{figure}

These are read directly from Fig. \ref{fig:seconddoubtB}, in a way similar to the previous pair
of links.

\section{A more general question: the $hgqi_u^d$-classes of 3-manifolds}
\begin{figure}[!h]
\begin{center}
\includegraphics[width=16cm]{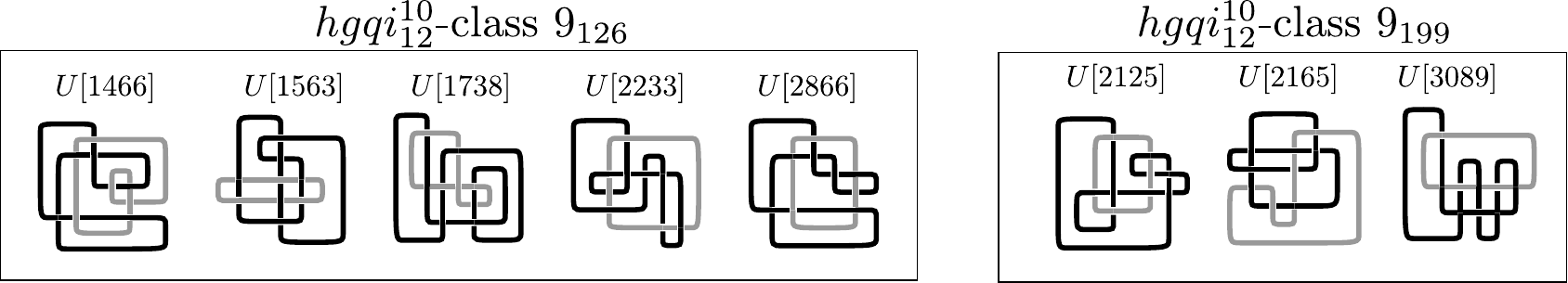}
\caption{\sf {\bf Note's final challenge:}  classify topologically $9_{126}$ and $9_{199}$. 
Here, to classify has the following strict meaning: for each pair 
of closed oriented hyperbolic 3-manifolds induced by links in one of these classes,
either make available a 
homeomorphism between them or, in the hyperbolic case, make available an isomophism between
their fundamental groups, or else make available an invariant which distinguishes them. Such
a proof of the coincidence or distinctiveness must be computationally short 
and reproducible by other researchers. The given projections
define blackboard framed links and so, the (integer) framing of each component is its self-writhe. 
In a blackboard framed
link the algorithm to get the presentation for the fundamental group of the associated 3-manifold
does not need and, thus, does not use the framing. Moreover, any integer framing can be realized
as a blackboard framed link by introducing appropriate curls in the projections to adjust the sef-writhes.
GAP \cite{gap2002gap} and SnapPea
(\cite{weeks2001snappea}) are good softwares to distinguish 
manifolds, but we personally have not tried them yet. 
It is a simple matter to obtain a canonical gem with $8n$ vertices from a blackboard framed link
with $n$ crossings, \cite{lins2007blink}. 
Gems are good at displaying homeomorphism via $TS$- and $u^n$-moves,
\cite{lins1995gca}. It factors the homeomorphism as a sequence of blob cancellations 
and valid flips (\cite{lins2006blobs}), never increasing the number of
vertices of the gems. Because of the lexicography 
inherent to graph with edges properly colored, 
a gem-based homeomorphism between two 3-manifolds will coincide
in any independent implementation of the algorithm given in Chapter 4 of \cite{lins1995gca}: 
the sequences of blob and flips turn out to be exactly the same.
If the manifolds are homeomorphic, of course each possible invariant will fail to distinguish
them. Therefore, to prove that two framed links are indeed manifestation of the same 
manifold we must make available a homeomorphism; or in the  hyperbolic case, to make available
an isomorphism between the fundamental
groups. To establish an explicit homeomorphism, 
what else could be used beyond a (short) path in a graph whose vertices are gems
and whose edges are either a blob cancellation or a valid flip (very simple moves)? 
Moreover, such an answer has the virtue of being quickly verifiable by independent implementations.
Is there a substitute for gems to acomplish 
this task? Kirby's moves \cite{kirby1978calculus} and their variants by Fenn and Rourke
\cite{fenn1979kirby} and more recently Martelli \cite{martelli2011finite}, are, with taylored 
exceptions, unusable because they increase the size of the links in completely blind directions and so, 
helplessly inferior to gems in this regard. The presentation of 3-manifolds 
based on the special spines of Matveev \cite{matveev2007algorithmic} 
seems to be a possibility, but first a theory to
deal with completion of the census and the isomorphism problem 
of such spines, as well as using some filter on them to decrease redundance,
is yet to be established and made available. In the case of gems the corresponding theory is simpler and 
is available since 1995, \cite{lins1995gca}.}
\label{fig:nine126nine199}
\end{center}
\end{figure}

The 3-manifolds of \cite{lins2007blink} are classified by homology and the quantum WRT$_r$-invariants
$r=3,\ldots,u$, up to $d$ decimal digits forming $hgqi_u^d$-classes. Our algorithm for computing the 
$WRT_r^d$-invariants are based on the theory developed in \cite{kauffman1994tlr}. 
The actual values rely on independent implementations
which coincide throughout \cite{kauffman1994tlr} and \cite{lins2007blink}.
The main domain 
of links in \cite{lins2007blink} (there are others) is formed by the so called {\em representative g-blinks},
$U[p]$'s $p=1,2,\ldots,$ which is a highly filtered class of blackboard framed links indexed by lexicography.
An important result of the work is that the $U[p]$'s form a universal class of 3-manifolds, in the sense 
that no closed orientable 3-manifold is missing. The examples of the previous 
section embed into two $hgqi_{12}$-classes:
$9_{126}$ (page 201 of  \cite{lins2007blink}) and $9_{199}$ (page 213 of \cite{lins2007blink}). 
The $hgqi_{12}^{10}$-class
$9_{126}$ is formed by 5 links $U[1466], U[1563], U[1738], U[2233]$ and $U[2866]$. 
The $hgqi_{12}^{10}$-class $9_{199}$ is formed by 3 links: $U[2125], U[2165]$ and $U[3089]$.
In Fig. \ref{fig:nine126nine199}, we display $9_{126}$ and $9_{199}$.
This note's final challenge is to classify topologically $9_{126}$ and $9_{199},$ in the sense
given in the caption of Fig. \ref {fig:nine126nine199}.

\section{Definition of Gem}
For completeness we briefly recall the basic definitions of gem theory, leading to its definition, \cite{lins1995gca}.
A {\em 4-graph} $G$ is a finite bipartite 4-regular graph whose edges are partitioned into 4 colors,
0,1,2, and 3, 
so that at each vertex there is an edge of each color, a proper edge-coloration, \cite{bondy1976gta}.
For each $i \in \{0,1,2,3\}$, let $E_i$ denote the set of $i$-colored edges of $G$.
A $\{j,k\}$-residue in a $4$-graph $G$ is a connected component of the subgraph induced by $E_j \cup E_k$.
A 2-residue is a $\{j,k\}$-residue, for some distinct colors $j$ and $k$.
A {\em gem} is a 4-graph $G$ such that for each color $i$, $G\backslash E_i$ can be embedded in the plane 
such that the boundary of each face is a 2-residue. From a gem there exists a straightforward
algorithm to obtain a closed orientable 3-manifold, in two different, dual ways. Every such a manifold is obtainable in this way.
An unecessary big gem is obtained from a triangulation $T$ for a manifold by taking the dual of the 
barycentric subdivision of $T$. Here the colors corresponds to the dimensions. Doing simplifications in the gem
completely destroys this correspondence.

\section{Conclusion}
A closed orientable 3-manifold is denoted {\em $n$-small} if it is induced by surgery on
a blackboard framed link with at most $n$ crossings.
Our bet is that both pairs of 3-manifolds in the 2 first sections of 
this short note are not homeomorphic. This would mean that the $9$-small manifolds are
completely classified and that
the combinatorial dynamics of Chapter 4 in \cite{lins1995gca} based 
on $TS$-moves which leads to a  (small, in the case of hyperbolic 3-manifolds)
number of minimal gems, named the {\em attractor of
the 3-manifold} is successful. This induces an efficient algorithm which 
is capable of classifying topologically all the 3-manifolds given as a blackboard framed link
with up to (so far) 9 crossings and maintains live the two Conjectures of page 15 of \cite{lins1995gca}:
the $TS$- and $u^n$-moves yield an efficient algorithm
to classify $n$-small 3-manifolds by explicitly displaying homeomorphisms, whenever they exist.

\bibliographystyle{plain}
\bibliography{bibtexIndex.bib}

\begin{thebibliography}{10}

\bibitem{bondy1976gta}
J.A. Bondy and U.S.R. Murty.
\newblock {\em {Graph theory with applications}}.
\newblock Macmillan London, 1976.

\bibitem{fenn1979kirby}
R.~Fenn and C.~Rourke.
\newblock On {K}irby's calculus of links.
\newblock {\em Topology}, 18(1):1--15, 1979.

\bibitem{gap2002gap}
GAP Group et~al.
\newblock Gap --- {G}roups, {A}lgorithms, and {P}rogramming, version 4.3, 2002.

\bibitem{kauffman1991knots}
L.H. Kauffman.
\newblock {\em Knots and physics}, volume~1.
\newblock World Scientific Publishing Company, 1991.

\bibitem{kauffman1994tlr}
L.H. Kauffman and S.~Lins.
\newblock {Temperley-Lieb Recoupling Theory and Invariants of 3-manifolds}.
\newblock {\em Annals of Mathematical Studies, Princeton University Press},
  134:1--296, 1994.

\bibitem{kirby1978calculus}
R.~Kirby.
\newblock A calculus for framed links in ${S}^3$.
\newblock {\em Inventiones Mathematicae}, 45(1):35--56, 1978.

\bibitem{lickorish1962representation}
W.B.R. Lickorish.
\newblock A representation of orientable combinatorial 3-manifolds.
\newblock {\em Annals of Mathematics}, 76(3):531--540, 1962.

\bibitem{lins2007blink}
L.D. Lins.
\newblock Blink: a language to view, recognize, classify and manipulate
  3{D}-spaces.
\newblock {\em Arxiv preprint math/0702057}, 2007.

\bibitem{lins1995gca}
S.~Lins.
\newblock {\em {Gems, Computers, and Attractors for 3-Manifolds}}.
\newblock World Scientific, 1995.

\bibitem{lins2006blobs}
S.~Lins and M.~Mulazzani.
\newblock {Blobs and flips on gems}.
\newblock {\em Journal of Knot Theory and its Ramifications}, 15(8):1001--1035,
  2006.

\bibitem{martelli2011finite}
B.~Martelli.
\newblock A finite set of local moves for {K}irby calculus.
\newblock {\em Arxiv preprint arXiv:1102.1288}, 2011.

\bibitem{matveev2007algorithmic}
S.~Matveev.
\newblock {\em Algorithmic topology and classification of 3-manifolds},
  volume~9.
\newblock Springer, 2007.

\bibitem{matveev2008}
S.~Matveev.
\newblock Personal communication during the workshop in {L}ow {D}imensional
  {T}opology, {O}berwolfach, {G}ermany, 2008.

\bibitem{rolfsen2003knots}
D.~Rolfsen.
\newblock {\em Knots and links}.
\newblock American Mathematical Society, 2003.

\bibitem{stillwell1993classical}
J.~Stillwell.
\newblock {\em {Classical Topology and Combinatorial Group Theory}}.
\newblock Springer Verlag, 1993.

\bibitem{weeks2001snappea}
J.~Weeks.
\newblock Snap{P}ea: a computer program for creating and studying hyperbolic
  3-manifolds, 2001.

\bibitem{whitehead1997}
A.~N. Whitehead.
\newblock {\em The classical group: their invariants and representations}.
\newblock Princeton {U}niversity {P}ress, 1997.

\end{thebibliography}

\vspace{5mm}
\begin{center}
\hspace{7mm}
\begin{tabular}{l}
   S\'ostenes L. Lins\\
   Centro de Inform\'atica, UFPE \\
   Av. Jornalista Anibal Fernandes s/n\\
   Recife, PE 50740-560 \\
   Brazil\\
   sostenes@cin.ufpe.br
\end{tabular}
\hspace{20mm}
\hspace{7mm}
\begin{tabular}{l}
Lauro D. Lins\\
AT\&T Labs Research \\
180 Park Avenue \\
Florham Park, NJ 07932 \\
USA\\
llins@research.att.com
\end{tabular}
\hspace{20mm}

\end{center}

\end{document}